\documentclass[11pt]{amsart}
\usepackage{amsmath,amssymb,amsfonts,url}

\numberwithin{equation}{section}

\newtheorem{thm}{Theorem}
\renewcommand{\thethm}{\kern -.8ex}
\newtheorem{lem}{Lemma}[section]

\newtheorem{prop}{Proposition}[section]
\newtheorem{cor}{Corollary}[section]

\begin{document}
\title{Prescribing curvatures on three dimensional
Riemannian manifolds with boundaries} 

\author{Lei Zhang}
\address{Department of Mathematics\\
        University of Alabama at Birmingham\\
        1300 University Blvd, 452 Campbell Hall\\
        Birmingham, Alabama 35294-1170}
\email{leizhang@math.uab.edu}
\thanks{Lei Zhang is supported by National Science Foundation Grant 0600275 (0810902)}

\date{\today}

%%%%%%%%%%%%%%%%%%%%%%%%%%%%%%%%%%%%%%%%%%%%%
\begin{abstract}
Let $(M,g)$ be a complete three dimensional Riemannian manifold
with boundary $\partial M$. Given smooth functions $K(x)>0$ and
$c(x)$ defined on $M$ and $\partial M$, respectively, it is
natural to ask whether there exist metrics conformal to $g$ so
that under these new metrics, $K$ is the scalar curvature and $c$
is the boundary mean curvature. All such metrics can be described
by a prescribing curvature equation with a boundary condition.
With suitable assumptions on $K$,$c$ and $(M,g)$ we show that all
the solutions of the equation can only blow up at finite points
over each compact subset of $\bar M$, some of them may appear on
$\partial M$. We describe the asymptotic behavior of the blowup
solutions around each blowup point and derive an energy estimate
as a consequence.
\end{abstract}
%%%%%%%%%%%%%%%%%%%%%%%%%%%%%%%%%%%%%%%%%%%%%

\maketitle

\section{Introduction}

In geometric analysis two well known problems are closely related,
the Nirenberg problem (or the Kazdan-Warner problem) and the
Yamabe problem. The Nirenberg problem asks what function $K(x)$ on
$\mathbb S^n$ is the scalar curvature of a metric $g$ on $\mathbb
S^n$ conformal to the standard metric $g_0$. From another point of
view, on any given compact Riemannian manifold $(M,g)$ without
boundary, the Yamabe problem, which was solved through the works
of Yamabe \cite{yamabe}, Trudinger\cite{trudinger},
Aubin\cite{aubin} and Schoen \cite{schoen}, concerns whether it is
possible to deform $g$ conformally to get a new metric with
constant scalar curvature. Similar questions can also be asked on
general Riemannian manifolds with boundaries: Let $(M,g)$ be a
Riemannian manifold with boundary, let $K\in C^1(M)$ and $c\in
C^1(\partial M)$ be $C^1$ functions defined on $M$ and $\partial
M$, respectively, then is it possible to deform $g$ conformally to
another metric $g_1$ so that $K$ is the scalar curvature and $c$
is the boundary mean curvature under $g_1$? By writing $g_1$ as
$g_1=u^{\frac 4{n-2}}g$, $K$ and $c$ are related to the scalar
curvature $R_g$ and the mean curvature $h_{g}$ under metric $g$ by
\begin{equation}\label{bry1}
\left\{\begin{array}{ll}
K=-\frac{4(n-1)}{n-2}u^{-\frac{n+2}{n-2}}(\Delta_gu-\frac{n-2}{4(n-1)}R_gu).
\\
c=\frac 2{n-2}u^{-\frac{n}{n-2}}(\partial_{
\nu_g}u+\frac{n-2}2h_gu).
\end{array}
\right.
\end{equation}
where $\nu_g$ is the unit normal vector pointing to the outside of
$\partial M$. $\Delta_g$ is the Laplace-Beltrami operator, which
can be written as $\Delta_g=\frac{1}{\sqrt g}\partial_i(\sqrt g
g^{ij}\partial_j)$ in local coordinates.

 If $K$ and $c$ are constants and $(M,g)$ is compact, the
existence of a metric with constant scalar curvature and boundary
mean curvature is always referred to as the boundary Yamabe
problem (BYP). Many cases of the BYP have been solved by Escobar
\cite{JE1,JE2,JE4}, Han-Li\cite{HL1} and Marques \cite{FM}. But
unlike the completely solved Yamabe problem, some cases of the BYP
are still open. In general, the boundary terms in (\ref{bry1})
make the nature of equation (\ref{bry1}) very different from its
counterpart without the boundary condition. Some difficulties
created by the boundary terms are still not completely understood.
In this article, we focus on three dimensional Riemannian
manifolds with boundaries and consider the corresponding
prescribing curvature equations defined on these manifolds.

 Let $(M,g)$ be a smooth three dimensional complete Riemannian manifold with boundary $\partial M$.
Suppose $K(x)>0$ and $c(x)$ are $C^1$ functions defined on $M$ and $\partial M$, respectively.
If $g_1=u^4g$ ($u>0$ smooth) is a metric conformal to $g$
that takes $K$ as the scalar curvature and $c$
as the mean curvature on $\partial M$, one can write the equation as
\begin{equation}
\label{aug31e1}
\left\{\begin{array}{ll}
\Delta_gu-\frac 18R_gu+\frac 18K(x)u^5=0,\quad M,\\
\\
\partial_{\nu_g}u+\frac 12h_gu=\frac 12c(x)u^3,\quad \partial M,
\end{array}
\right.
\end{equation}
The purpose of this article is to understand the bubbling
phenomena of (\ref{aug31e1}) under natural assumptions on $K$ and
$c$. Without the boundary condition, the bubbling phenomenon of
equation (\ref{aug31e1}) and its high dimensional variants have
been studied extensively under various assumptions. The reader is
referred to the following incomplete list and the references
therein
\cite{berti,SB,SB1,SB2,CSL1,CSL2,CSL3,Druet1,Druet2,Ah2,Khuri,lipart1,LZhang,LZ1,LZ2,LZhu,FM,S1,schoen,Z3}.
However, much less work can be found to address the case with
boundary conditions. In this article we describe the blowup
phenomenon for (\ref{aug31e1}) under weak assumptions on $K$ and
$c$.

To state our main result, we first remark that $M$ may not be compact, so our
result concerns the blowup phenomenon for every compact subset of $M$ that shares with
$M$ a part of its boundary. Namely, let $V$ be a subset of $M$ such
that $\bar V$ is compact and let $\Gamma=\partial V\cap \partial M$ be the part of the
boundary that $V$ shares with $M$, then our main result can be stated as follows:
\begin{thm}
\label{thm5}
Let $(M,g)$, $R_g$, $h_g$, $V$, $\Gamma$, $K$, $c$ be described as above.
Let $u>0$ be a classical solution of (\ref{aug31e1}). Suppose $\Gamma$ is umbilic and there exists
$\epsilon_0>0$ such that
for every point $p\in \Gamma$, $g$ is conformal to the Euclidean metric $\delta$ in
$B(p,\epsilon_0)\cap M$. Assume in addition that for some $\Lambda>0$,
$$\frac 1{\Lambda}\le K(x)\le \Lambda,\quad  M$$
$$\max\bigg(\|\nabla K(x)\|_{L^{\infty}(M)},\|c\|_{C^1(\partial M)}\bigg )\le \Lambda.$$
Then for some $C(M,g,\Lambda,V,\epsilon_0)>0$ and integer $m(M,g,\Lambda,V,\epsilon_0)>0$,
there exist
local maximum points of $u$, denoted as $S:=\{P_1,...,P_m\}$, such that
$$\mbox{dist}_g(P_i,P_j)\ge \frac 1C,\quad \frac 1Cu(P_i)\le u(P_j)\le Cu(P_i),
\quad \forall i\neq j,$$
\begin{equation}
\label{sep7e1}
\frac 1C\sum_{l=1}^m\xi_{P_l,u(P_l)^2}\le u\le C\sum_{l=1}^m\xi_{P_l,u(P_l)^2},
\quad \mbox{on}\quad \bar V
\end{equation}
where
$$\xi_{Q,\mu}(P):=(\frac{\mu}{1+\mu^2\mbox{dist}(P,Q)^2})^{\frac 12},\quad P\in \bar V.$$
\end{thm}

As a consequence, we have the following energy estimate:
\begin{cor}
\label{thm6}
Under the assumptions of Theorem \ref{thm5}, we have
\begin{equation}
\label{aug31e3}
\int_V(|\nabla_gu|^2+u^6)\le C(M,g,\Lambda,V,\epsilon_0).
\end{equation}
\end{cor}

If $(M,g)$ is compact, Theorem \ref{thm5} concludes that if the metric is locally conformally
flat near $\partial M$, which is umbilic, then there are only finite blowup
point in $\bar M$ and the solutions of (\ref{aug31e1}) can be estimated by (\ref{sep7e1}).

The results of Theorem \ref{thm5} and the Corollary \ref{thm6} are closely related to a priori
estimates of the solutions. When $(M,g)$ is compact with a boundary and $K$ and $c$ are constants,
various compactness results and related discussions
can be found in \cite{Ara}, \cite{Ah2}, \cite{HL1}, \cite{Ah} and the reference therein.
The $C^1$ assumptions on $K$ and $c$ in Theorem \ref{thm5} should be sharp. The reason
is even for the interior equations the $C^1$ assumption on $K$
is necessary to have a description of the blowup phenomenon by (\ref{sep7e1}). So Theorem
\ref{thm5} implies that for this three dimensional situation, the $C^1$ assumption on
$c$ does not lead to more restrictions on $K$. For dimensions $4$
and higher, we expect the situation to be much
more subtle, first the assumption on $K$ will be more delicate and should be made
in neighborhoods of its critical points, then the sign of $c$ at the
blowup point will make a more significant difference. Also, the flatness assumptions in
the neighborhood of the critical points of $K$ and $c$ might be related. As
the dimension grows higher, the relations will become more subtle, new phenomenon will
come out.

The major step in the proof of \ref{thm5} is to establish a
Harnack type inequality near $\Gamma$. This is done by the well
known method of moving planes. Unlike most of the previous works
on this type of Harnack inequality established only for the
interior equations, the equation in Theorem \ref{thm5} has new
features that were not dealt with before. Generally speaking, the
method of moving planes requires a delicate construction of some
test functions. Now this becomes more difficult because of two
reasons, first, the symmetry of the domain is destroyed, second,
when blowup happens on the boundary or close to the boundary, the
function $c$ creates more error terms for the test functions to
control. In this article we found a way to handle all the
difficulties for dimension three, our approach is motivated by
Caffarelli's magnificent ideas in \cite{ca1}\cite{ca2} on free
boundary problems.

 In the next few sections, we mainly focus on the proof of the Harnack type inequality on
upper half Euclidean balls and we shall mention how to obtain the interior estimate away
from the boundary $\Gamma$.

\section{Harnack type inequalities on the boundary}

The proof of Theorem \ref{thm5} can be divided into two parts. First we consider the region
close to $\Gamma$: $N(\Gamma, \epsilon_0)$. i.e. the $\epsilon_0$ neighborhood of $\Gamma$.
In this region we shall establish a Harnack type inequality. Then outside this region, we
use the techniques developed in \cite{LZ1} to get the same type of Harnack inequality. Note that
for this interior part we don't need to assume the metric to be locally conformally flat.

Since $(M,g)$ is locally conformally flat near $\Gamma$, at each
point $p\in \Gamma$ we can find a positive smooth function $\phi$ such that
$\delta=\phi^4g$ where $\delta$ is the Euclidean metric. The umbilicity of $\Gamma$ implies that
$\partial M$ near $p$ is either a piece of sphere or a hyperplane in the new metric. Since
both neighborhoods are conformally equivalent through the inversion map $f(x)=x|x|^{-2}$, we
just assume that near $p$ the boundary is a hyperplane, then
the equation in the neighborhood of $p$ can be written as
$$\left\{\begin{array}{ll}
\Delta (\frac{u}{\phi})+K(\frac{u}{\phi})^5=0&\qquad B^+_{\epsilon_1}\subset \mathbb R^3 \\ \\
\partial_{x_3}(\frac{u}{\phi})=-c(\frac{u}{\phi})^3
&\qquad \mbox{on}\quad \partial 'B^+_{\epsilon_1}.
\end{array}
\right.
$$
where $\epsilon_1>0$ is a positive constant that only depends on
$\epsilon_0$, $M,g,\Gamma$. $B_{\epsilon_1}^+$ is the upper half
ball of radius $\epsilon_1$, $\partial'B_{\epsilon_1}^+=\partial
B_{\epsilon_1}^+\cap \partial \mathbb R^3_+$. The computation
above is based on the following conformal covariant properties of
the two operators:
\begin{eqnarray*}
\Delta (\frac{u}{\phi})&=&\phi^{-5}(\Delta_{g}u-c(n)R_gu)\\
-\partial_{x_3}(\frac{u}{\phi})&=&\phi^{-3}(\partial_{\nu_g}u+\frac 12h_gu).
\end{eqnarray*}

So the whole equation is reduced to the Euclidean case. From now on in this section
we just consider the following case (by abusing the notations we still use $K$
to denote a positive $C^1$ function on a upper half disk $B^+_{3R}$ and $c$
to denote a $C^1$ function defined $\partial'B^+_{3R}$, the lower part of
$\partial B^+_{3R}$):

\begin{equation}
\label{chap4e1}
\left\{\begin{array}{ll}
\Delta u+K(x)u^5=0&\qquad B^+_{3R}\subset \mathbb R^3 \\ \\
\displaystyle{\frac{\partial u}{\partial x_3}}=c(x')u^3
&\qquad \mbox{on}\quad \partial 'B^+_{3R}.
\end{array}
\right.
\end{equation}
where $\partial 'B^+_{3R}=\partial B^+_{3R}\cap \partial \mathbb
R^3_+$, $x'$ is the projection of $x$ to $\partial 'B^+_{3R}$. The
main result of this section is
\begin{prop}
\label{thm1}
Let $u\in C^1(\overline{B^+_{3R}})\cap C^2(B^+_{3R})$ be a positive function that solves
(\ref{chap4e1}). Assume for some $\Lambda>0$, $K$ and $c$ satisfy
\begin{equation}
\label{KK}
\frac 1{\Lambda}\le K(x)\le \Lambda,\quad
|\nabla K(x)|\le \Lambda  R^{-1}, \,\, \forall x\in \overline{B^+_{3R}}
\end{equation}
and
\begin{equation}
\label{cc}
|c(x')|\le \Lambda,\quad |\nabla c(x')|\le \Lambda R^{-1},\quad
\forall x'\in \partial'B^+_{3R}
\end{equation}
respectively, then for some $C(\Lambda)>0$
\begin{equation}
\label{mar23e1}
\max_{\bar B_R^+}u\cdot \min_{\bar B_{2R}^+}u \le C(\Lambda)R^{-1}.
\end{equation}
\end{prop}

This Harnack type inequality reveals important information about the interaction of
bubbles (large local maximum points) of $u$ which can be seen in the following

\begin{cor}
\label{mar30cor1}
Under the assumptions in Theorem \ref{thm1}, there exist a positive constant $C(\Lambda)$,
a positive integer $m(\Lambda)$, and a set of $m$ local maximum points of $u$ in $\bar
B^+_R$, denoted as $S:=\{P_1,..,P_m\}$, such that
$$\mbox{dist}(P_i,P_j)\ge R/C,\quad \frac 1Cu(P_j)\le u(P_i)\le Cu(P_j),\quad
\forall i\neq j,$$
\begin{equation}
\label{sep12e1}
\frac 1C\sum_{l=1}^m\xi_{P_l,u(P_l)^2}\le u\le C\sum_{l=1}^m\xi_{P_l,u(P_l)^2},
\quad \mbox{on}\quad \bar B^+_R
\end{equation}
where
$$\xi_{Q,\mu}(P):=(\frac{\mu}{1+\mu^2\mbox{dist}(P,Q)^2})^{\frac 12},\quad P\in
\bar B^+_R.$$
Consequently $u$ satisfies
\begin{equation}
\label{sep12e2}
\int_{B_R^+}(|\nabla u|^2+u^6) \le C(\Lambda).
\end{equation}
\end{cor}

In the proof of Proposition we shall omit a selection process of Schoen since it is well known
to experts. We only include it in the appendix. At the end of this section we
indicate the outline of the proof of Corollary \ref{mar30cor1} based on Proposition \ref{thm1}.
The proof can be found in \cite{LZhang} with obvious changes, but we indicate
the key points for the convenience of the readers.

\subsection{The Proof of Proposition \ref{thm1}}
We only need to consider the case $R=1$. The general case can be reduced to the case $R=1$
by considering the function $R^{\frac 12}u(R\cdot)$.
The proof is by a contradiction. Suppose (\ref{mar23e1})
does not hold, there exists a sequence $u_i(x_i)$ such that
\begin{equation}
\label{mar25e1}
u_i(x_i)\min_{\bar B^+_2}u_i\ge i
\end{equation}
where $u_i(x_i)=\max_{\bar B^+_1}u_i$ and $x_i\in \bar B_1$. Clearly the above inequality implies
$u_i(x_i)\to \infty$.
By a standard selection process of
Schoen \cite{S1} and the classification theorems of Caffarelli-Gidas-Spruck \cite{CGS}
and Li-Zhu \cite{LZhu}, one can consider $x_i$ as
a local maximum of $u_i$, moreover the following sequence of functions
\begin{equation}
\label{apr2e1}
\bar v_i(y):=u_i(x_i)^{-1}u_i(u_i(x_i)^{-2}y+x_i)
\end{equation}
 converge in $C^2$ over any finite domain
in the following two cases:
\begin{enumerate}
\item If $\lim_{i\to \infty}u_i^2(x_i)x_{i3}\to \infty$, $\bar
v_i(y)$ converge in $C^2_{loc}(\mathbb R^3)$ to $U(y)$ which
satisfies
\begin{equation}
\label{apr3e1} \left\{\begin{array}{ll} \Delta U(y)+\lim_{i\to
\infty}K_i(x_i)U(y)^5=0,\quad \mathbb R^3,
\\
U(0)=1=\max_{\mathbb R^3}U.
\end{array}
\right.
\end{equation}
\item
If $\lim_{i\to \infty}u_i^2(x_i)x_{i3}$ is bounded, by passing to a subsequence we assume
$$\lim_{i\to \infty}u_i^2(x_i)x_{i3}=t_0.$$
In this case $\bar v_i$ converge uniformly to $U$ which satisfies
\begin{equation}
\label{apr2e2}
\left\{\begin{array}{ll}
\Delta U(y)+\lim_{i\to \infty}K_i(x_i)U(y)^5=0,&\quad y_3>t_0,\\
\partial_3U(y)=\lim_{i\to \infty}c_i(x_i')U^3,&\quad y_3=-t_0,\\
U(0)=1=\max U.
\end{array}
\right.
\end{equation}
\end{enumerate}
Since the selection process and the application of the classification theorems are standard.
We put the details in the Appendix. Similar techniques can also be found in
\cite{CSL1}, \cite{LZhang}, \cite{Z3}, etc.

The proof of Proposition \ref{thm1} that follows can be divided into three steps. First we
rule out case one. i.e. We shall show that the blowup points can not be
far away from the boundary.
In the second step we prove the case of $\lim_{i\to \infty}c(x_i')\le 0$. Then in the final step
we prove the case of $\lim_{i\to \infty}c(x_i')>0$.

\bigskip

\subsubsection{Step one} In this subsection we derive a contradiction to
$$\lim_{i\to \infty}T_i:=\lim_{i\to \infty}u^2_i(x_i)x_{i3}\to \infty.$$
With no loss of generality we assume $\lim_{i\to \infty}K_i(x_i)=3$, in this case
the $U$ in (\ref{apr3e1}) is of the form
$$U(y)=(1+|y|^2)^{-\frac 12}.$$ Let $M_i=u_i(x_i)$. In this subsection we use $v_i$ in stead of
$\bar v_i$. i.e.
$$v_i(y)=M_i^{-1}u_i(M_i^{-2}y+x_i).$$
The selection process implies $B(x_i,\frac 12)\cap \mathbb
R^3_+\subset B^+_2$ so the rescaled domain for $v_i$ will have a
part of the boundary (the upper part) whose distance to $0$ is
comparable to $M_i^2$. By (\ref{mar25e1}) we know on this part of
the boundary $v_i(y)|y|>i$. By choosing $\epsilon_i\to 0$ slowly
we can still obtain $v_i(y)|y|>\sqrt i$ for $|y|=\epsilon_iM_i^2$.
So we define
$$\Omega_i=\{y:\,M_i^{-2}y+x_i\in B(x_i,\epsilon_i)\cap \mathbb R^3_+\} $$
and
$$
\partial'\Omega_i=\partial \Omega_i\cap \{y_3=-T_i\},\quad
\partial''\Omega_i=\partial \Omega_i\setminus \partial'\Omega_i.$$
Then $v_i$ satisfies
\begin{equation}
\left\{\begin{array}{ll}
\Delta v_i(y)+K_i(M_i^{-2}y+x_i)v_i(y)^5=0
\qquad y\in \Omega _i \\
\\
\partial_3 v_i(y)
\le \Lambda v_i^3\qquad y\in \partial \Omega _i\cap \{y_3=-T_i\},\quad T_i:=M_i^2x_{i3}
\end{array}
\right.
\label{f6a1}
\end{equation}
Note that in this case we don't need the specific equation for $v_i$ on $\partial '\Omega _i$,

By the discussion on the behavior of $v_i$ on $\partial''\Omega_i$ we have
\begin{equation}
\label{mar25e2} v_i(y)|y|\to \infty \qquad \mbox{for}\quad y\in
\partial''\Omega_i.
\end{equation}
Let $\lambda\in [\frac 12, 2]$ and
$$
{v_i}^{\lambda}(y):=(\frac{\lambda}{|y|})v_i({\frac{\lambda^2y}{|y|^2}}),
$$ then
direct computation gives
\begin{equation}
\label{f6z1}
 \Delta {v_i}^{\lambda}(y)+K_i(M_i^{-2}y^{\lambda}+x_i)
v_i^{\lambda}(y)^5=0 \qquad \mbox{for}\ \ |y| >\lambda/2.
\end{equation}
Here we use $y^{\lambda}=\lambda^2y/|y|^2$.

In this step we assume $\lambda\in [\frac 12, 2]$. Let $\Sigma_{\lambda}=\Omega_i\setminus
\bar B_{\lambda}$. Note that for simplicity we shall omit $i$ in some notations.
Let $w_{\lambda}=v_i-v_i^{\lambda}$ and we consider the equation for $w_{\lambda}$:
\begin{equation}
\left\{\begin{array}{ll}
\Delta w_{\lambda}+b_{\lambda}w_{\lambda}=Q_{\lambda},&\quad \Sigma_{\lambda},\\ \\
\partial_3w_{\lambda}\le \Lambda \xi_iw_{\lambda}+\Lambda (v_i^{\lambda})^3-\partial_3(v_i^{\lambda}),
&\quad \partial'\Omega_i\setminus \bar B_{\lambda}.
\end{array}
\right.
\label{mar25e3}
\end{equation}
where $b_{\lambda}$ and $\xi_i$ are obtained from mean value theorem:
\begin{eqnarray*}
b_{\lambda}(y)&=&5K_i(M_i^{-2}y+x_i)\int_0^1(tv_i(y)+(1-t)v_i^{\lambda}(y))^4dt \\
\xi_i(y)&=&3\int_0^1(tv_i+(1-t)v_i^{\lambda})^2dt.
\end{eqnarray*}
and
$$Q_{\lambda}(y)=(K_i(M_i^{-2}y^{\lambda}+x_i)-K_i(M_i^{-2}y+x_i))(v_i^{\lambda}(y))^5$$
is estimated as follows:
\begin{equation}
|Q_{\lambda}(y)|\le C_0(\Lambda)M_i^{-2}|y|^{-4},\quad \Sigma_{\lambda}.
\label{mar25e4}
\end{equation}

Since $v_i$ converges in $C^2$ norm over $U$ over any fixed finite
domain, $v_i^{\lambda}$ is close to $U^{\lambda}$, the Kelvin
transformation of $U$. By direct computation $U>U^{\lambda}$ for
$|y|>\lambda$ and $\lambda \in (0,1)$. On the other hand
$U<U^{\lambda}$ for $\lambda>1$ and $|y|>\lambda$. So the strategy
of the proof is to find a test function $h_{\lambda}$ ($i$ omitted
in this notation) so that the moving sphere method works for
$w_{\lambda}+h_{\lambda}$, and the $h_{\lambda}$ is just a
perturbation of $w_{\lambda}$, which means
$h_{\lambda}(y)=\circ(1)|y|^{-1}$ in $\Sigma_{\lambda}$. Then it
is possible to move the spheres from a position less than $1$ to a
position larger than $1$ keeping $w_{\lambda}+h_{\lambda}>0$ in
$\Sigma_{\lambda}$. But this is a contradiction since
$w_{\lambda}+h_{\lambda}$ converges to $U-U^{\lambda}$ in finite
domains and $U<U^{\lambda}$ for $\lambda>1$ and $|y|>\lambda$.

The test function in this section is
\begin{equation}
\label{mar25e5}
h_{\lambda}(y)=C_0M_i^{-2}(-\frac 12|y|^{-2}+\lambda^{-1}|y|^{-1}-\frac 12\lambda^{-2}).
\end{equation}
$h_{\lambda}$ is a radial function, the function $h_{\lambda}(r)$ satisfies
$$\left\{\begin{array}{ll}
h_{\lambda}''(r)+\frac{2}rh'_{\lambda}(r)=-C_0M_i^{-2}r^{-4},\quad r>\lambda, \\
h_{\lambda}(\lambda)=h_{\lambda}'(\lambda)=0.
\end{array}
\right.
$$
By maximum principle, $h_{\lambda}(r)<0$ for $r>\lambda$. Note that the $C_0$ in (\ref{mar25e5})
is the same as the one in (\ref{mar25e4}).

Now we consider the equation for $w_{\lambda}+h_{\lambda}$, from (\ref{mar25e3}) we have
\begin{equation}
\label{mar30e3}
(\Delta +b_{\lambda})(w_{\lambda}+h_{\lambda})=Q_{\lambda}+\Delta h_{\lambda}+b_{\lambda}h_{\lambda}
\le 0\quad \Sigma_{\lambda}.
\end{equation}
The reason is $\Delta h_{\lambda}+Q_{\lambda}\le 0$ and $h_{\lambda}\le 0$ in $\Sigma_{\lambda}$.

On $\partial'\Omega_i\setminus \bar B_{\lambda}$ we have
\begin{equation}
\label{mar25e6}
\partial_3(w_{\lambda}+h_{\lambda})\le \Lambda\xi_i(w_{\lambda}+h_{\lambda})+\Lambda(v_i^{\lambda})^3
-\partial_3(v_i^{\lambda})+\partial_3h_{\lambda}-\Lambda\xi_ih_{\lambda}.
\end{equation}

By the definition of $h_{\lambda}$ one can verify that $h_{\lambda}=\circ(1)|y|^{-1}$ over
$\Sigma_{\lambda}$. Let
$$O_{\lambda}=\{y\in \Sigma_{\lambda};\quad v_i(y)\le 2v_i^{\lambda}(y)\,\, \}.$$ Then we observe
that the maximum principle in the moving sphere process only needs to be applied over $O_{\lambda}$
because outside $O_{\lambda}$, $w_{\lambda}+h_{\lambda}>0$.
Now we claim that
\begin{equation}
\label{mar25e7}
\Lambda(v_i^{\lambda})^3-\partial_3(v_i^{\lambda})+\partial_3h_{\lambda}-\Lambda\xi_ih_{\lambda}<0,
\quad \mbox{on}\quad
\partial' \Omega_i\cap \overline{O_{\lambda}}.
\end{equation}

Once we have (\ref{mar25e7}), (\ref{mar25e6}) becomes
\begin{equation}
\label{mar30e2}
\partial_3(w_{\lambda}+h_{\lambda})\le \Lambda\xi_i(w_{\lambda}+h_{\lambda})
\quad \mbox{on}\quad \partial'\Omega_i\cap \overline{O_{\lambda}},
\end{equation}
which is the form for the application of the moving sphere method.

To see (\ref{mar25e7}), we first observe that
\begin{eqnarray*}
\partial_3(v_i^{\lambda})&=&\lambda |y|^{-3}T_iv_i(y^{\lambda})+
2\lambda^3|y|^{-5}\sum_{j=1}^3\partial_jv_i(y^{\lambda})y_jT_i\\
&&+\lambda^3|y|^{-3}\partial_3v_i(y^{\lambda})
\quad \mbox{on}\quad \partial'\Omega_i\setminus \bar B_{\lambda}.
\end{eqnarray*}
The second term and the third term on the right are negligible
comparing to the first term on the right. The reason is $v_i\to U$
in $C^2_{loc}(\mathbb R^3)$, which implies
$$\nabla v_i(y^{\lambda})\to \nabla U(0)=0.$$
On the other hand, by assumption $T_i\to \infty$, we
see the first term dominates the other two terms.

Next we see that the assumption $T_i\to \infty$ makes $|\partial_3(v_i^{\lambda})|$ dominate
$\Lambda(v_i^{\lambda})^3$ as the later is of the order $O(|y|^{-3})$. By the definition of $h_{\lambda}$
and the definition of the domain $\Omega_i$, we see easily that $|\partial_3(v_i^{\lambda})|$
also dominates $\partial_3h_{\lambda}$. Finally in $\bar O_{\lambda}$, $\xi_i$ is comparable to
$|y|^{-2}$, so
$\xi_ih_{\lambda}$ is comparable to $O(M_i^{-2}|y|^{-2})$, which is much smaller than
$|\partial_3(v_i^{\lambda})|$, (\ref{mar25e7}) is proved.

The process of making the moving sphere process start is standard, even though the boundary
condition makes it different from the interior case.

\begin{lem}
\label{mar26lem1}
For any fixed $\lambda_0\in (\frac 12, 1)$ and all large $i$,
$$w_{\lambda_0}+h_{\lambda_0}>0 \quad \mbox{in}\quad \Sigma_{\lambda_0}.$$
\end{lem}

\noindent{\bf Proof of Lemma \ref{mar26lem1}:} Since $v_i\to U$ in
$C^2_{loc}(\mathbb R^3)$ and $U>U^{\lambda}$ for $|y|>\lambda$ if
$\lambda\in (\frac 12, 1)$. By the convergence of $v_i$ to $U$ and
the fact $h_{\lambda}=\circ(1)$ over finite domains we can check
easily that for any fixed $R_1>>1$,
$$w_{\lambda_0}+h_{\lambda_0}>0\quad \mbox{in}\quad \Sigma_{\lambda_0}\cap B_{R_1}.$$
Since $R_1$ is chosen sufficiently
large, one can find $\sigma (\lambda_0)>0$ such that
$$v_i^{\lambda_0}(y)\le (1-3\sigma)|y|^{-1} \quad |y|\ge R_1$$
and
$$v_i(y)>(1-\sigma)|y|^{-1}\quad \mbox{on}\quad |y|=R_1.$$
The definitions of $h_{\lambda_0}$ and $\Omega_i$ imply $|h_{\lambda_0}(y)|\le
\frac{\sigma}5|y|^{-1}$ for $|y|>R_1$. So to finish the proof of this Lemma
it is enough to show
\begin{equation}
\label{mar26e1}
v_i(y)>(1-2\sigma)|y|^{-1}\quad \mbox{in}\quad \Omega_i\setminus \overline{B_{R_1}}.
\end{equation}

On $|y|=R_1$ and $\partial''\Omega_i$ we certainly have $v_i>(1-2\sigma)|y|^{-1}$.
To apply the maximum principle to the super harmonic function
$v_i-(1-2\sigma)|y|^{-1}$, we need
\begin{equation}
\label{apr10e1}
\partial_3(v_i-(1-2\sigma)|y|^{-1})\le \Lambda (v_i^3-(1-2\sigma)^3|y|^{-3}).
\quad \partial'\Omega_i
\end{equation}
Using $T_i\to \infty$ we obtain by direct computation that
$$\partial_3((1-2\sigma)|y|^{-1})>\Lambda((1-2\sigma)|y|^{-1})^3
\quad \partial'\Omega_i.$$
(\ref{apr10e1}) follows immediately.
(\ref{apr10e1}) means the non-positive minimum can not be attained on $\partial'\Omega_i$.
Lemma \ref{mar26lem1} is established. $\Box$

\bigskip

Lemma \ref{mar26lem1} means the moving sphere process can start at $\lambda=\lambda_0$. The purpose
of the moving spheres is to show that
\begin{equation}
\label{mar30e1}
w_{\lambda}+h_{\lambda}>0\,\, \mbox{in}\,\, \Sigma_{\lambda}\quad \mbox{for all }
\lambda\in [\lambda_0,\lambda_1].
\end{equation}
Once this is proved, we have $U>U^{\lambda_1}$ for some
$|y|>\lambda_1$, which is a contradiction to the fact that $\lambda_1>1$.
To see why (\ref{mar30e1}) holds, first, (\ref{mar30e3}) and (\ref{mar30e2}) means
the maximum principle holds for $w_{\lambda}+h_{\lambda}$, second, (\ref{mar25e2}) means there is no touch
on $\partial''\Omega_i$. Step one is established. $\Box$

\subsubsection{Step two}

In this step, we deal with the case $\lim_{i\to \infty}c_i(x_i')\le 0$. Let
$$-c_0=\lim_{i\to \infty}c_i(x_i')\le 0.$$
In this case the limit of $\bar v_i$ is $U$, which satisfies
$$\left\{\begin{array}{ll}
\Delta U+3U^5=0 \quad y\in \mathbb R^3, y_3>-\lim_{i\to \infty}T_i, \\
\partial_3U=-c_0U^3 \quad y_3=-\lim_{i\to \infty}T_i,\\
U(0)=1=\max U.
\end{array}
\right.
$$
By Li-Zhu's classification theorem \cite{LZhu}, $\lim_{i\to \infty}T_i=0$.
So we define $v_i$ as
$$v_i(y)=M_i^{-1}u_i(M_i^{-2}y+x_i')$$
where $x_i'$ is the projection of $x_i$ on $\partial \mathbb
R^3_+$. Then $v_i\to U$ in $C^2_{loc}(\bar{\mathbb R}^3_{+})$.
Since $U(0)=1$, we know, from Li-Zhu's theorem, that
\begin{equation}
\label{mar27e1}
U(y)=(1+c_0^2)^{\frac 12}(1+(1+c_0^2)^2|y'|^2+(1+c_0^2)(y_3+\frac{c_0}{1+c_0^2})^2)^{-\frac 12}.
\end{equation}
Let $v_i^{\lambda}(y)=(\frac{\lambda}{|y|})v_i(y^{\lambda})$. In this step we let
$$\lambda_0=\frac 12(1+c_0^2)^{-\frac 12} ,\quad \lambda_1=2(1+c_0^2)^{-\frac 12}$$ and we
require
$\lambda\in [\lambda_0, \lambda_1]$ and the moving sphere method will be applied on
$\Sigma_{\lambda}=B^+(0,\epsilon_iM_i^2)\setminus \bar B_{\lambda}$.
Direct computation gives
\begin{equation}
\label{mar27e2}
\left\{\begin{array}{ll}
\Delta v_i^{\lambda}(y)+K_i(M_i^{-2}y^{\lambda}+x_i')(v_i^{\lambda}(y))^5=0,\quad \Sigma_{\lambda},\\
\\
\partial_3v_i^{\lambda}(y)=c_i(M_i^{-2}y^{\lambda}+x_i')(v_i^{\lambda}(y))^3,\quad
\partial \Sigma_{\lambda}\cap \partial \mathbb R^3_+
\end{array}
\right.
\end{equation}
Also by direct computation one can verify that
$$\left\{\begin{array}{ll}
U(y)>U^{\lambda}(y),&\quad |y|>\lambda,\quad 0<\lambda<(1+c_0^2)^{-\frac 12},\\
\\
U(y)<U^{\lambda}(y),&\quad |y|>\lambda,\quad \lambda>(1+c_0^2)^{-\frac 12}
\end{array}
\right.
$$
For this reason we require $\lambda\in [\lambda_0,\lambda_1]$ in this step.
The equation for $w_{\lambda}=v_i-v_i^{\lambda}$ becomes
\begin{equation}
\label{mar27e3}
\left\{\begin{array}{ll}
\Delta w_{\lambda}+K_i(M_i^{-2}\cdot +x_i')\xi_1w_{\lambda}=Q_1,\quad \Sigma_{\lambda},\\
\\
\partial_3w_{\lambda}=c_i(M_i^{-2}\cdot +x_i')\xi_2w_{\lambda}+Q_2,\quad \partial \Sigma_{\lambda}
\cap \partial \mathbb R^3_+.
\end{array}
\right.
\end{equation}
where $\xi_1$ and $\xi_2$ are obtained by mean value theorem:
$$\xi_1=5\int_0^1(tv_i+(1-t)v_i^{\lambda})^4dt,\quad
\xi_2=3\int_0^1(tv_i+(1-t)v_i^{\lambda})^2dt$$
and $Q_1$ and $Q_2$ are error terms to be controlled by test functions:
\begin{eqnarray*}
Q_1(y)&=&(K_i(M_i^{-2}y^{\lambda}+x_i')-K_i(M_i^{-2}y+x_i'))v_i^{\lambda}(y)^5,\quad \Sigma_{\lambda},\\
Q_2(y)&=&(c_i(M_i^{-2}y+x_i')-c_i(M_i^{-2}y^{\lambda}+x_i'))v_i^{\lambda}(y)^3,\quad
\partial \Sigma_{\lambda}\cap \partial \mathbb R^3_+.
\end{eqnarray*}
For $Q_1$ and $Q_2$ we use
\begin{eqnarray}
|Q_1(y)| &\le & C_2(\Lambda)M_i^{-2}|y|^{-4},\quad \Sigma_{\lambda}, \nonumber \\
|Q_2(y)| &\le & C_2(\Lambda )M_i^{-2}|y|^{-2}(1-\frac{\lambda}{|y|}),\quad
\partial \Sigma_{\lambda}\cap \partial \mathbb R^3_+.
\label{mar27e4}
\end{eqnarray}
Note that in $Q_2$ it is important to have the $1-\frac{\lambda}{|y|}$ term, even though
we don't need this term to appear in the estimate of $Q_1$.

The construction of the test function here consists of two parts, $h_1$ and $h_2$.
We first define
\begin{equation}
\label{mar28e5}
h_1(r)=C_3(\Lambda)(r^{-\frac 12}-\frac 12\lambda^{\frac 12}r^{-1}-\frac 12\lambda^{-\frac 12})
M_i^{-2}
\end{equation}
where $C_3(\Lambda)>>C_2$ is to be determined. $h_1(r)$ is the solution of
$$\left\{\begin{array}{ll}
h_1''(r)+\frac 2rh_1'(r)=-\frac{C_3}4M_i^{-2}r^{-\frac 52},\quad r>\lambda,\\
h_1(\lambda)=h_1'(\lambda)=0.
\end{array}
\right.
$$
By maximum principle, $h_1(r)<0$ for $r>\lambda$. $h_1(|y|)$ is the first part of the test
function $h_{\lambda}$. To define the second part of the test function, we let $\phi:[1,\infty)\to
[0,\infty)$ be a
smooth non-negative function that satisfies
$$\left\{\begin{array}{ll}
\phi(1)=0,\quad \phi(r)=r^{-2}\quad \mbox{for}\quad r\ge 2,\\
\phi'(1)>0, \quad \phi(r)>0\quad \mbox{for}\quad r>1.
\end{array}
\right.
$$
Let $\phi_{\lambda}(\cdot)=\phi(\lambda \cdot)$ and we define
\begin{equation}
\label{mar28e6}
h_2(y)=-C_4(\Lambda)M_i^{-2}y_3\phi_{\lambda}(|y|)
\end{equation}
where $C_4(\Lambda)$ is chosen so large that
\begin{equation}
\label{mar28e1}
C_4\phi_{\lambda}(|y|)>2C_2|y|^{-2}(1-\frac{\lambda}{|y|}),\quad |y|>\lambda.
\end{equation}
Note that on $\partial'\Sigma_{\lambda}:=\partial
\Sigma_{\lambda}\cap \partial \mathbb R^3_+$,
\begin{equation}
\label{mar28e2}
\partial_3h_2(y)=-C_4M_i^{-2}\phi_{\lambda}(|y|), \quad \partial_3h_1(y)=0 \quad
\partial'\Sigma_{\lambda}.
\end{equation}
Let $h_{\lambda}=h_1+h_2$, the equation for $w_{\lambda}+h_{\lambda}$ is
\begin{eqnarray*}
(\Delta+K_i(M_i^{-2}\cdot+x_i')\xi_1)(w_{\lambda}+h_{\lambda})
=Q_1+\Delta h_{\lambda}+K_i(M_i^{-2}\cdot+x_i')\xi_1h_{\lambda}\quad \Sigma_{\lambda}, \\
(\partial_3-c_i(M_i^{-2}\cdot +x_i')\xi_2)(w_{\lambda}+h_{\lambda})=
-c_i(M_i^{-2}\cdot +x_i')\xi_2h_{\lambda}+\partial_3h_2+Q_2\quad \partial'\Sigma_{\lambda}.
\end{eqnarray*}
For the application of the moving sphere method we show that the right hand sides of the above
are non-positive. Namely we shall show
\begin{equation}
\label{mar28e3}
Q_1+\Delta h_1+\Delta h_2+K_i(M_i^{-2}\cdot+x_i')\xi_1h_{\lambda}\le 0 \quad \Sigma_{\lambda}
\end{equation}
and
\begin{equation}
\label{mar28e4}
-c_i(M_i^{-2}\cdot +x_i')\xi_2h_{\lambda}+\partial_3h_2+Q_2\le 0 \quad \partial'\Sigma_{\lambda}
\cap \overline{O_{\lambda}}
\end{equation}
where $O_{\lambda}=\{y\in \Sigma_{\lambda};\,\, v_i(y)\le 2v_i^{\lambda}(y)\,\, \}$ is the only
place where the maximum principle needs to hold, because by the definition of $h_{\lambda}$,
$h_{\lambda}(y)=\circ(1)|y|^{-1}$ in $\Sigma_{\lambda}$, which means $w_{\lambda}+h_{\lambda}>0$ in
$\Sigma_{\lambda}\setminus \overline{O_{\lambda}}.$

To see (\ref{mar28e3}), first by (\ref{mar27e4}) and (\ref{mar28e5}) one sees that if $C_3$ is
large enough
$$\Delta h_1+Q_1\le -\frac 18 C_3M_i^{-2}|y|^{-2.5},\quad \Sigma_{\lambda}.$$
For $h_2$ we have
$$\Delta h_2=-C_4M_i^{-2}(\phi_{\lambda}''(|y|)+\frac 4{|y|}\phi_{\lambda}'(|y|))y_3
\quad \Sigma_{\lambda}.$$
So by choosing $C_3$ larger if necessary we have
$$\Delta (h_1+h_2)+Q_1\le -\frac 1{16}C_3M_i^{-2}|y|^{-2.5} \quad \Sigma_{\lambda}.$$
The term $K_i(M_i^{-2}\cdot +x_i)\xi_1h_{\lambda}$ is non-positive. So (\ref{mar28e3}) is proved.

To see (\ref{mar28e4}), first by the definition of $h_2$ we have
$$\partial_3h_2+Q_2<-\frac{C_4}2M_i^{-2}\phi_{\lambda}(|y|), \quad \partial'\Sigma_{\lambda}.$$
Since $\lim_{i\to \infty}c_i(x_i')=-c_0\le 0$, one uses the fact $|y|\le \epsilon_iM_i^2$ and the
uniform bound of the $C_1$  norm of $c_i$ to get
\begin{eqnarray*}
&&-c_i(M_i^{-2}\cdot+x_i')\xi_2h_{\lambda}\\
&\le &\circ(1)\xi_2|h_1+h_2| \le \circ(1) M_i^{-2}
|y|^{-2}(1-\frac{\lambda}{|y|})\quad \partial'\Sigma_{\lambda}\cap \bar O_{\lambda}.
\end{eqnarray*}
The second inequality above is because in $\bar O_{\lambda}$, $\xi_2\sim |y|^{-2}$. We also
use $h_{\lambda}(\lambda)=0$ in the above. Then (\ref{mar28e4}) follows immediately.

Next we show that the moving sphere process can get started, namely we have
\begin{lem}
\label{mar28lem1}
$w_{\lambda_0}+h_{\lambda_0}>0$ in $\Sigma_{\lambda_0}$ for all large $i$.
\end{lem}

\noindent{\bf Proof of Lemma \ref{mar28lem1}:}
The proof is similar to the one in step one. For $|y|\le R_1$ ($R_1>>1$), $w_{\lambda_0}
+h_{\lambda_0}>0$ is guaranteed by the expressions of $U$, $U^{\lambda_0}$ and the
convergence of $v_i$ to $U$ over finite domains. Also $h_{\lambda_0}=\circ(1)|y|^{-1}$ means
$h_{\lambda_0}$ is just a perturbation over finite domains. Moreover on $|y|=R_1$, there
is $\epsilon_0>0$ such that
$$U(y)-U^{\lambda_0}(y)>5\epsilon_0|y-e_3|^{-1},\quad |y|=R_1,\quad e_3=(0,0,1).$$
Also we have
$$v_i^{\lambda_0}(y)\le (1-4\epsilon_0)(1+c_0^2)^{-\frac 12}
|y-e_3|^{-1},\quad |y|>R_1, y\in \Sigma_{\lambda_0}$$
and
$$|U^{\lambda_0}(y)-v_i^{\lambda_0}(y)|\le \epsilon_0|y-e_3|^{-1},\quad |y|>R_1$$
for all $i$ large.
So for $y\in \Sigma_{\lambda}\setminus B_{R_1}$ we only need to show
$$v_i(y)>(1-\epsilon_0)(1+c_0^2)^{-\frac 12}|y-e_3|^{-1}\quad y\in \Sigma_{\lambda}\setminus B_{R_1}.$$
Since $v_i(y)-(1-\epsilon_0)(1+c_0^2)^{-\frac 12}|y-e_3|^{-1}$ is super-harmonic function in
$\Sigma_{\lambda}\setminus B_{R_1}$ and is positive on $|y|=R_1$ and $|y|=\epsilon_iM_i^2$. The only
place to consider is $\partial'\Sigma_{\lambda_0}\setminus B_{R_1}$, on which we have
\begin{eqnarray}
\label{apr10e2}
&&\partial_3(v_i(y)-(1-\epsilon_0)(1+c_0^2)^{-\frac 12}|y-e_3|^{-1})\nonumber \\
&=&c_i(M_i^{-2}y+x_i')v_i^3(y)
-(1-\epsilon_0)(1+c_0^2)^{-\frac 12}|y-e_3|^{-3}.
\end{eqnarray}
Since $\lim_{i\to \infty}c_i(x_i')\le 0$, (\ref{apr10e2}) means
$$\partial_3(v_i(y)-(1-\epsilon_0)(1+c_0^2)^{-\frac 12}|y-e_3|^{-1})<0$$
at places where $v_i$ is close to $(1-\epsilon_0)(1+c_0^2)^{-\frac 12}|y-e_3|^{-1}$
on $\partial'\Sigma_{\lambda_0}\setminus B_{R_1}$. This means
$$v_i>(1-\epsilon_0)(1+c_0^2)^{-\frac 12}|y-e_3|^{-1}\quad \mbox{on}\quad
\partial'\Sigma_{\lambda_0}\setminus B_{R_1}.$$
Lemma \ref{mar28lem1} is established. $\Box$

Once we have Lemma \ref{mar28lem1}, the moving sphere process can start at $\lambda=\lambda_0$.
Note that the equation for $w_{\lambda}+h_{\lambda}$ becomes
$$\left\{\begin{array}{ll}
(\Delta+K_i(M_i^{-2}\cdot+x_i')\xi)(w_{\lambda}+h_{\lambda})\le 0,\quad \Sigma_{\lambda},
\\ \\
(\partial_3-c_i(M_i^{-2}\cdot +x_i')\xi)(w_{\lambda}+h_{\lambda})\le 0,\quad
\partial'\Sigma_{\lambda}\cap \bar O_{\lambda}.
\end{array}
\right.
$$
This means the maximum principle always holds for $w_{\lambda}+h_{\lambda}$ as long as it is
positive on the boundary $\partial''\Omega_i$, which is certainly the case because $v_i>>v_i^{\lambda}$
on $\partial'\Omega_i$ and $h_{\lambda}$ is a perturbation. So the spheres can be moved to
$\lambda_1=2(1+c_0^2)^{-\frac 12}$ to get a contradiction.

\subsubsection{Step three}
In this subsection we deal with the case $\lim_{i\to
\infty}c_i(x_i')=c_0>0$. $\bar v_i$ converges uniformly to $U$ in
all finite subsets of $\{y\in \mathbb R^3;\,\, y_3\ge -\lim_{i\to
\infty}M_i^2x_{i3}\,\, \}$. $U$ satisfies (still assuming
$\lim_{i\to \infty}K_i(x_i)=3$)
$$\left\{\begin{array}{ll}
\Delta U+3U^5=0,\quad y_3>-\lim_{i\to \infty}M_i^2x_{i3} \\
\partial_3U=c_0U^3,\quad y_3=-\lim_{i\to \infty}M_i^2x_{i3},\\
U(0)=1.
\end{array}
\right.
$$
By Li-Zhu's classification theorem $\lim_{i\to \infty}M_i^2x_{i3}=c_0$, and
$$U(y)=(1+|y|^2)^{-\frac 12}.$$
Let
$$v_i(y)=M_i^{-1}u_i(M_i^{-2}y+x_i').$$
Then
$$v_i\to U_1(y)=(1+|y-c_0e_3|^2)^{-\frac 12}.$$
Let $v_i^{\lambda}$ be the Kelvin transformation of $v_i$. In this case we let
$$\lambda_0=\frac 12(1+c_0^2)^{\frac 12},\quad \lambda_1=2(1+c_0^2)^{\frac 12}$$
because
\begin{eqnarray*}
U_1>U_1^{\lambda}\quad \mbox{for}\quad \lambda\in (0,(1+c_0^2)^{\frac 12}),\quad |y|>\lambda,\\
U_1<U_1^{\lambda}\quad \mbox{for}\quad \lambda>(1+c_0^2)^{\frac 12}),\quad |y|>\lambda.
\end{eqnarray*}
Let $w_{\lambda}=v_i-v_i^{\lambda}$, then the equation for $w_{\lambda}$ is still
described by (\ref{mar27e3}) with $\xi_1$, $\xi_2$, $Q_1$,$Q_2$ defined as in step two and the
estimates for $Q_1$ and $Q_2$ are still (\ref{mar27e4}). Since $\lim_{i\to \infty}
c_i(x_i')=c_0>0$, the construction of the test function $h_{\lambda}$ is this
case is much more delicate.

The construction of $h_{\lambda}$ consists of two parts. First we shall construct $h_1$
to control the region close to $\partial B_{\lambda}$, then this function is connected
smoothly to $0$. Clearly $h_1$ creates new difficulties at the regions where it is connected
to $0$. Next we use $h_2$ to control the region far enough to $\partial B_{\lambda}$.
$h_2$ is $0$ in regions close to $\partial B_{\lambda}$ and becomes negative as $y$ is
far away from $\partial B_{\lambda}$. A parameter of $h_2$ is chosen to be large so that
$h_2$ not only controls the difficulties from $Q_1$ and $Q_2$, but also those
from $h_1$. One delicate thing is that $h_2$ does not create new difficulties. By
choosing all the parameters carefully we shall obtain the following properties for $h_{\lambda}$:

\begin{eqnarray}
&&\Delta h_{\lambda}+K_i(M_i^{-2}y+x_i')\xi_1h_{\lambda}+Q_1\le 0\quad \Sigma_{\lambda}
\nonumber \\
&&\partial_3h_{\lambda}-c_i(M_i^{-2}y+x_i')\xi_2h_{\lambda}+Q_2\le
0, \quad \partial \Sigma_{\lambda}\cap \partial \mathbb R^3_+\cap
\overline{O_{\lambda}}.
\nonumber \\
&&h_{\lambda}=0\quad \mbox{on}\quad \partial B_{\lambda}\cap
\mathbb R^3_+,
\nonumber \\
&&h_{\lambda}=\circ(1)|y|^{-1}\quad \Sigma_{\lambda}
\label{apr27e1}
\end{eqnarray}

We first construct $h_2$. Let
$$\phi(y)=1-A^{\frac 14}|y+re_3|^{-\frac 14},\quad r=|y|$$
where $A>5\lambda_1$ is to be determined. Let $\Omega_1\subset
\mathbb R^3_+$ be the region where $\phi$ is positive. Let
$\Gamma_1$ be the $0$ level surface of $\phi$. Direct computation
shows that
$$\mathbb R^3_+\setminus B(0,2A)\subset \Omega_1\subset \mathbb R^3_+\setminus B(0,A/3).$$
Let
$$h_2(y)=\left\{\begin{array}{ll}
0&\qquad \mathbb R^3_+\setminus (\Omega_1\cup B_{\lambda})\\
-C_6M_i^{-2}\phi^2&\qquad \Sigma_{\lambda}\cap \Omega_1
\end{array}
\right.
$$
where $C_6$ is a large constant to be determined. From here we see $h_2$ is $C^1$ across
$\Gamma_1$. By direct computation
\begin{eqnarray}
\Delta h_2&=&-C_6M_i^{-2}(2|\nabla \phi|^2+2\phi \Delta \phi)\nonumber \\
&=&-C_6M_i^{-2}\bigg (
\frac 18A^{\frac 12}(2r^2+2ry_3)^{-\frac 94}(5r^2+8ry_3+3y_3^2) \nonumber \\
&&+\frac 12\phi A^{\frac 14}(2r^2+2ry_3)^{-\frac {17}8}(\frac 34r^2+2ry_3+\frac 54y_3^2)\bigg )
\nonumber \\
&\le &-\frac{C_6}{50}M_i^{-2}A^{\frac 12}r^{-\frac 52} \quad \Omega_1\cap \Sigma_{\lambda}.
\label{apr27e2}
\end{eqnarray}
From here we observe that $h_2$ is weakly super-harmonic in
$\Sigma_{\lambda}$ since $\nabla h_2=0$ on $\Gamma_1$. On
$\partial \Omega_1\cap \partial \mathbb R^3_+\cap
\overline{O_{\lambda}}$ we need to estimate
$$\partial_3h_2-c_i(M_i^{-2}\cdot +x_i')\xi_2h_2+Q_2.$$ Since
$$c_i(M_i^{-2}y+x_i')\xi_2\le \frac{M}{1+r^2}\quad \mbox{on}\quad \overline{O_{\lambda}}
\cap \partial \mathbb R^3_+$$ for some $M>0$ depending on $c_0$,
it is enough to show
\begin{equation}
\label{apr27e3}
\partial_3h_2-\frac{M}{1+r^2}h_2+Q_2\le 0,\quad
\partial \Omega_1\cap \partial \mathbb R^3_+\cap \overline{O_{\lambda}}.
\end{equation}

By the definition of $h_2$, if $A(c_0)$ is large, we have
\begin{eqnarray}
&&\partial_3h_2-\frac{M}{1+r^2}h_2=-2C_6M_i^{-2}(\partial_3\phi-\frac{M}{2(1+r^2)}\phi)\phi
\nonumber \\
&\le &-2C_6M_i^{-2}(\partial_3\phi-\frac{M}{2(1+r^2)})\phi \nonumber\\
&=&-2C_6M_i^{-2}(\frac 14A^{\frac 14}(2r^2)^{-\frac 98}r-\frac{M}{2(1+r^2)})\phi\nonumber \\
&\le &-\frac{C_6}8M_i^{-2}A^{\frac 14}r^{-\frac 54}\phi \quad
\partial \Omega_1\cap \partial \mathbb R^3_+. \label{apr27e5}
\end{eqnarray}
Here $A$ is fixed. We observe that if
$|y|>16A$, $\phi(y)>\frac 12$. So by choosing $C_6$ large enough we have
\begin{equation}
\label{apr27e4}
(\partial_3-c_i(M_i^{-2}\cdot+x_i')\xi_2)h_2+Q_2\le
-\frac{C_6}{16}M_i^{-2}A^{\frac 14}r^{-\frac 54} \quad \partial
(\Sigma_{\lambda}\setminus B_{16A})\cap \partial \mathbb R^3_+\cap
\overline{O_{\lambda}}.
\end{equation}
By (\ref{apr27e5}) we also know
\begin{equation}
\label{apr27e6} (\partial_3-c_i(M_i^{-2}\cdot+x_i')\xi_2)h_2\le 0
\quad \partial (\Sigma_{\lambda}\cap B_{16A})\cap \partial \mathbb
R^3_+\cap \overline{O_{\lambda}}.
\end{equation}
So the most important feature of $h_2$ is it does not create new difficulties, even
though it does not control all the error terms in the whole $\Sigma_{\lambda}$.

Now we define $h_1$ to be
\begin{equation}
\label{apr27e7}
h_1(y)=\left\{\begin{array}{ll}
-C_5M_i^{-2}(r-\lambda)e^{Ny_3}&\quad \Sigma_{\lambda}\cap B_{16A}\\
0&\quad \Sigma_{\lambda} \setminus B_{17A},\\
\mbox{smooth connection such that $h_1\le 0$ in }&\quad
\Sigma_{\lambda}\cap \overline{B_{17A}}\setminus B_{16A}
\end{array}
\right.
\end{equation}
where $C_5$ and $N$ are large positive constants to be determined. $A$ is the one in
the definition of $h_2$, which has been determined. Also the "smooth connection" in
the definition means $|\nabla h_1(y)|,|\nabla^2h_1(y)|\le C_7(C_5,A)M_i^{-2}$
for $16A<|y|<17A$, $y\in \Sigma_{\lambda}$.

By the definition of $h_1$,
\begin{eqnarray}
\Delta h_1(y)&=&-C_5M_i^{-2}e^{Ny_3}(\frac 2r+\frac {2N}ry_3+N^2(r-\lambda))
\nonumber \\
&\le &-2C_5M_i^{-2}r^{-1}\quad y\in B_{16A}\cap \Sigma_{\lambda}.
\label{apr27e8}
\end{eqnarray}
Here we require $C_5$ to satisfy
\begin{equation}
\label{apr27e9}
\Delta h_1+Q_1\le 0\quad y\in B_{16A}\cap \Sigma_{\lambda}.
\end{equation}
On the boundary we have
\begin{eqnarray*}
&&\partial_3h_1-c_i(M_i^{-2}y'+x_i')\xi_2h_1\\
&=&-C_5M_i^{-2}(r-\lambda)(N-c_i(M_i^{-2}y'+x_i')\xi_2) \quad
\partial (\Sigma_{\lambda}\cap B_{16A})\cap \partial \mathbb
R^3_+\cap \overline{O_{\lambda}}.
\end{eqnarray*}
So by choosing $N(c_0)$ large and $C_5$ larger if necessary we have
\begin{equation}
\partial_3h_1-c_i(M_i^{-2}y'+x_i')\xi_2h_1+Q_2\le 0\quad
\partial (\Sigma_{\lambda}\cap B_{16A})\cap \partial \mathbb R^3_+\cap\overline{O_{\lambda}}.
\end{equation}

So far $C_5,N,A$ are all determined, we finally determine $C_6$. One last requirement
for $C_6$ is to control the bad part of $h_1$ in $(B_{17A}\setminus B_{16A})\cap \Sigma_{\lambda}$.
By choosing $C_6$ larger if necessary we see from (\ref{apr27e5}) and (\ref{apr27e4}) that
the errors caused by bending $h_1$ to $0$ can be controlled by $h_2$. Therefore
the first two equations of (\ref{apr27e1})
are established. Here we also use the fact that $h_{\lambda}\le 0$ in $\Sigma_{\lambda}$.
 By the definitions of $h_1$ and $h_2$, the other two equations of (\ref{apr27e1})
are also satisfied.
We are left with how to let the moving sphere process start.  This part is similar to
step two. We leave the details to the interested readers. Proposition \ref{thm1} is established.
$\Box$

\bigskip

Proposition \ref{thm1} is the major step in deriving the energy
estimate in Corollary \ref{mar30cor1}. Based on the previous works
of Y. Y. Li \cite{lipart1}, Han-Li \cite{HL1} and Li-Zhang
\cite{LZhang}, it is standard to derive Corollary \ref{mar30cor1}.
Therefore we only mention the major steps and the main idea in
this argument: First one uses Schoen's selection process to find
all large maximums of $u$ in, say, $B^+_{2R}$. Around each of
these local maximums there is a small neighborhood in which $u$
looks like a standard bubble which has most of its energy in it.
The distance between these local maximums is the crucial
information to find and this is the place where the Harnack
inequality in Proposition \ref{thm1} is used. The essential
difference between this locally defined equation and those
globally defined equations (such as the ones in \cite{lipart1}) is
that one can not find two bubbles closest to each other. For each
bubble, there certainly exists a bubble closest to it, but one
certainly can not assume the first bubble is the closest one to
the second bubble. This difficulty, which comes from the local
nature of the equation, requires a different approach than those
in \cite{lipart1}. The way to overcome this difficulty is to
rescale the equation so that after the scaling, the equation is
centered at the first bubble and the distance between the first
bubble and the second bubble is one. Then the Harnack inequality
in Proposition \ref{thm1} applied to the local region implies that
the first bubble and the second bubble must have comparable
magnitudes. Then it is possible to show that two bubbles can not
tend to the same blowup point because otherwise a harmonic
function with a positive second order term can be found. This
second term will lead to a contradiction in the Pohozaev Identity.
Note that the important second order term can only be proved to be
positive if the adjacent bubbles have comparable magnitudes, which
is the key information revealed by the Harnack inequality. Once we
have known that all bubbles are far apart, it is possible to use
standard elliptic estimates to show that the behavior of $u$ near
each large local maximum is like a harmonic function with fast
decay. So (\ref{sep12e1}) as well as (\ref{sep12e2}) can be
obtained. $\Box$

\subsection{Harnack inequality in the interior of $M$}
Now we consider $V\setminus N(\Gamma,\epsilon_0)$. Let $d$ be the distance
between $V\setminus N(\Gamma,\epsilon_0)$ and $\partial M$. Over this interior region we don't need
to assume $M$ to be locally conformally flat. We shall establish the following inequality:
\begin{prop}
\label{sep11p1}
Let $p\in V\setminus N(\Gamma,\epsilon_0)$, there exist $\epsilon_2(M,g,\Lambda,d)>0$ and
$C(M,g,\Lambda,d)>0$ such that
\begin{equation}
\sup_{B(p,\epsilon_2/3)}u\cdot \inf_{B(p,2\epsilon_2/3)}u\le C\epsilon_2^{-1}.
\label{sep11e1}
\end{equation}
\end{prop}

The proof of Proposition \ref{sep11p1} follows from the argument
in \cite{LZ1}. Only small modification is needed to adjust to the
current situation. The outline is as follows: Suppose
(\ref{sep11e1}) does not hold, then we can select $x_i$ as a
sequence of local maximums of $u_i$. Then the equation can be
written in a conformal normal coordinates centered at $x_i$. After
rescaling, $u_i$ becomes $v_i$, a sequence that converges in
$C^2_{loc}(\mathbb R^3)$ to a standard bubble whose maximum is
$1$. Then consider the Kelvin transformation of $v_i(y)$:
$v_i^{\lambda}(y)$. By comparing the equation for $v_i$ and
$v_i^{\lambda}$ we see that the only new term in this context is
the term:
$$(K_i(M_i^{-2}y^{\lambda}+x_i)-K_i(M_i^{-2}y+x_i))(v_i^{\lambda})^5 \quad \mbox{in}
\quad \Sigma_{\lambda}$$
where $\Sigma_{\lambda}$ becomes a symmetric domain defined appropriately. The above term
is of the harmless order $O(M_i^{-2}|y|^{-4})$. By using the same test function we used in
\cite{LZ1}, a contradiction can be obtained correspondingly. We leave the details to
the interested readers. Proposition \ref{sep11p1} is proved. $\Box$

\bigskip

Proposition \ref{thm1} and Proposition \ref{sep11p1} lead to Theorem \ref{thm5} by the
arguments in \cite{LZhang} and \cite{LZ2}.

\section{Appendix}
In this section we provide some details of the Schoen's selection process for the convenience
of readers. We only provide details for case one and case two. The details for case three
are similar.

\subsection{A Calculus Lemma}
First we use a Calculus lemma to simplify some computations. This Calculus lemma was first
used in \cite{LZhang}.
\begin{lem}
\label{f6lem4}
Let $u\in C^0(\overline{B_1\cap \{t>-T\}}) \ \  T\ge 0$.
Then for every $a>0$, there exists $x\in B_1\cap\{t\ge -T \}$ such that
$$u(x)\ge \frac{1}{2^a}\max _{\bar B^T_{\sigma}(x)}u, \quad \mbox{and}\quad
\sigma ^a u(x)\ge \frac{1}{2^a}u(0). $$
where $B^T_{\sigma}(x)=B(x,\sigma )\cap \{t\ge -T \}$, $\sigma =(1-|x|)/2 $.
\end{lem}

\noindent{\bf Proof of Lemma \ref{f6lem4}:}\ First we remark that if $T\ge 1$, the
selection is over the whole $B_1$. To prove the lemma, consider
$$v(y)=(1-|y|)^au(y). $$
Let $x\in \bar B^T_1 $ be a maximum point of $v$ and let $\sigma =
(1-|x|)/2$. By comparing $v(x)$ with $v(0)$ and $v(z)$ for all $z$ in $\bar B_{\sigma}^T(x)$,
we see that $x$ and $\sigma $ have the desired properties. $\Box$

\subsection{The selection process for case one in the proof of Proposition \ref{thm1}:}
In this subsection we explain why $x_i$ can be considered as a local maximum of $u_i$,
assuming $u_i^2(x_i)x_{i3}\to \infty$.
We shall apply Lemma \ref{f6lem4} for $T>1$. i.e. the selection is over
the whole $B_1$.
Let $r_i=|x_{i3}|/2$, so $u_i(x_i)r_i^{\frac 12}\to \infty$. Apply Lemma \ref{f6lem4}
to the function $u_i(x_i+r_i\cdot )$ over $B_1$, then the conclusion can be translated
as follows: there is $a_i\in B(x_i,r_i)$ such that
$$u_i(a_i)\ge \max_{\bar B(a_i,\sigma _i)} 2^{-\frac{1}{2}}u_i(x) $$
where $\sigma _i:=\frac{1}{2}(r_i-|a_i-x_i|)$ and
$$\sigma _i^{\frac{1}{2}}u_i(a_i)\ge 2^{-\frac{1}{2}}u_i(x_i)
r_i^{\frac{1}{2}} \to \infty $$
Let $$\tilde{v}_i(y)=u_i(a_i)^{-1}u_i(u_i(a_i)^{-2}y+a_i) \quad
|y|\le u_i(a_i)^{2}\sigma_i $$
Then $$\tilde{v}_i(y)\le 2^{\frac{1}{2}}\quad
\mbox{for}\ \ |y|\le u_i(a_i)^{2}\sigma_i\to \infty $$
Clearly $\tilde{v}_i$ satisfies
$$ \Delta \tilde{v}_i(y)+K_i(u_i(a_i)^{-2}y+a_i)
\tilde{v}_i(y)^{5}=0 \qquad |y| \le u_i(a_i)^{2}\sigma _i$$ So by
standard elliptic theory we know there is a subsequence of
$\tilde{v}_i(y)$ (still denoted by $\tilde{v}_i(y)$) that
converges uniformly to  $\tilde{U}_0(y)$ on all compact subsets of
$\mathbb R^3$. $\tilde{U}_0(y)$ satisfies

\begin{equation}
\Delta \tilde{U}_0(y)+\tilde{K}_0{\tilde{U}_0(y)}^{5}=0 \qquad
y\in \mathbb R^3 \label{Lin5}
\end{equation}
where $\tilde{K}_0=\lim_{i \to \infty} K_i(a_i)$. With no loss of generality
we assume $\tilde{K}_0=3$, so by the well known classification theorem
of Caffarelli-Gidas-Spruck \cite{CGS},
$$\tilde{v}_i(y)=(\frac{\mu}{1+\mu ^2|y-z|^2})^{\frac{1}{2}} $$
for some $z\in \mathbb R^3$ and $\mu \ge 1 $. So $\tilde{U}_0(y)$
has an absolute maximum at $z$. Consequently
$\tilde{v}_i(y)$ has a local maximum at $z_i$ close to $z$ when $i$ is large. \\
Let  $\bar{x}_i=u_i(a_i)^{-2}z_i+a_i$, then $\{\bar{x}_i\}$ are
local maximum points of $u_i$, also it is easy to verify that
$$u_i(\bar{x}_i)\ge \frac{\mu}{2}u_i(a_i)\ge \frac{\mu}{2}u(x_i) $$
and $\bar{x}_i\in B(a_i,\sigma _i)\subset B(x_i,r_i)$ which means
$\bar{x}_{i3}\ge \frac{1}{2}x_{i3}$, consequently
$$u_i^2(\bar{x}_i)\bar{x}_{i3}\to \infty $$
So we can consider $x_i$ as $\bar{x}_i$ at the beginning. The
$\bar v_i(y)$ defined in (\ref{apr2e1}) converges in
$C^2_{loc}(\mathbb R^3)$ to
$$U(y)=(1+|y|^2)^{-\frac 12}.$$

\subsection{Selection for case two in the proof of Proposition \ref{thm1}:}

In this subsection we show that $x_i$ can be considered as a local maximum of $u_i$ under
the assumption $u_i^2(x_i)x_{i3}=O(1)$.

By applying Lemma \ref{f6lem4} to $u_i(x_i+\frac 12\cdot )$ with $a=1/2$
and $T=-x_{i3}$ we can find
$\bar{x}_i\in \overline{B(x_i,1/2)}\cap \overline{\{x_3>-x_{i3}\}}$ such that
$$u_i(\bar{x}_i)\ge 2^{-\frac 12}\max u_i \qquad
\mbox{over}\quad \overline{B(x_i,\sigma_i)\cap \{x_3>-x_{i3}\} }$$
and
$\sigma_i^{\frac{1}{2}}u_i(\bar{x}_i)\ge
2^{-\frac{1}{2}}u_i(x_i)$
where $\sigma _i=\frac 12(\frac 12-|x_i-\bar{x}_i|)$. By these two inequalities we have
$$u_i(\bar{x}_i)\ge u_i(x_i)\to \infty ,\qquad
u_i(\bar{x}_i)\sigma_i^{\frac{1}{2}}\to \infty $$ and
$u_i(\bar{x}_i)\inf_{\partial B^+_{2}}u_i\ge i.$
Now we can assume that
$\bar{T}_i:=u_i^2(\bar{x}_i)\bar{x}_{i3}$ is bounded because otherwise by replacing $x_i$ by
$\bar x_i$ we go back to the case one, which has been discussed.
Let $$
\hat {v}_i(y)=u_i(\bar{x}_i)^{-1}u_i(u_i(\bar{x}_i)^{-2}y+\bar x_i).$$
Then $\hat {v}_i$ satisfies
$$\left\{\begin{array}{ll}
\Delta \hat {v}_i(y)+K_i(u_i(\bar{x}_i)^{-2}y+\bar{x}_i)
\hat v_i(y)^{5}=0 &\quad y\in \bar{\Omega}_i \\
\partial_3\hat {v}_i(y)=
c_i(u_i(\bar{x}_i)^{-2}y'+\bar{x}_i')\hat v_i(y)^3
&\quad \{y_3=-\bar{T}_i\}\cap \partial \bar{\Omega}_i
\end{array}
\right. $$ where $\bar{\Omega}_i=\{y;\quad
u_i(\bar{x}_i)^{-2}y+\bar{x}_i\in B^+_{3}\}$. It is clear that the
lower part of $\partial \bar{\Omega}_i$ is a subset of
$\{y_3=-\bar{T}_i\}$. Moreover by the facts
$u_i(\bar{x}_i)\bar{x}_{i3}^{\frac 12}\le C$ and
$u_i(\bar{x}_i)\sigma _i^{\frac{1}{2}}\to \infty $ we know $\hat
{v}_i$ is uniformly bounded on all compact subsets of $\{y_3\ge
-\lim_{i\to \infty}\bar{T}_i\}$. With no loss of generality we
assume that $K_i(\bar{x}_i)\to 3$, then by elliptic estimates,
$\hat v_i(y-\bar{T}_ie_3)$ converges in $C^2$ norm on all compact
subsets of $\bar{\mathbb R}^3_+$ to $\bar{U}$, which solves
$$\left\{\begin{array}{ll} \Delta \bar{U}+3\bar{U}^{5}
=0 &\qquad \mathbb R^3_+ \\
\partial_3 \bar{U}=c\bar{U}^3
&\qquad y_3=0
\end{array}
\right.
$$
where $c=\lim_{i\to \infty}c_i(\bar{x}_i')$. By Li-Zhu's
classification result \cite{LZhu}, we know that $\bar{U}$ has a
unique maximum point $z_0$ in $\overline{\mathbb R^3_+}$. By the
$C^2$ convergence of $\hat v_i(\cdot -\bar{T}_ie_3)$ to $\bar{U}$,
we can find $\{y_i\}_{i=1,2,..}$ as local maximum points of $\hat
{v}_i(\cdot -\bar{T}_ie_3)$ that approach $z_0$ as $i\to \infty$.
Then by the definition of $\hat v_i$, we know that $u_i(\bar
x_i)^{-2}(y_i-\bar{T}_ie_3)+\bar{x}_i$ are local maximum points of
$u_i$. So we redefine $x_i$ as
$u_i(\bar{x}_i)^{-2}(y_i-\bar{T}_ie_3)+\bar{x}_i$. So the $\bar
v_i$ defined in (\ref{apr2e1}) converges to the function $U_1$ in
(\ref{apr2e2}).

\end{document}